\RequirePackage{snapshot}
\documentclass{IEEEtran}
\usepackage{amsfonts,amsmath,amsthm,amssymb,dsfont}

\usepackage[usenames,dvipsnames]{color}

\ifCLASSINFOpdf 
\usepackage[pdftex]{graphicx}  
\usepackage[pdftex,colorlinks,
           bookmarks=true,%
           pdftitle=SINR\ in\ wireless\ networks\ and\ the\ two-parameter
  Poisson-Dirichlet\ process, 
           pdfauthor= H.\ P.\ Keeler\ -\ B.\ Blaszczyszyn%
           ]{hyperref}  
\usepackage[numbers,sort&compress]{natbib} 
\usepackage{hypernat} 
\else

 \usepackage[dvips]{graphicx} 
\usepackage[dvips, colorlinks,
           bookmarks=true,%
           pdftitle=SINR\ in\ wireless\ networks\ and\ the\ two-parameter
  Poisson-Dirichlet\ process, 
           pdfauthor= H. P.\ Keeler\ -\ B.\ Blaszczyszyn%
           ]{hyperref}
 \usepackage[numbers,sort&compress]{natbib} 
\fi 
\hypersetup{
   bookmarksnumbered,
   pdfstartview={FitH},
   citecolor={blue},
   linkcolor={red},
   urlcolor=[rgb]{0,0.55,0},
   pdfpagemode={UseOutlines}} 

\graphicspath{{./}{./Plots/}}

\newcommand{\E}{\mathbf{E}}

\newcommand{\calI}{\mathcal{I}}
\newcommand{\calJ}{\mathcal{J}}

\newcommand{\Prob}{\mathbf{P}}

\newcommand{\Ind}{\mathds{1}}
\newcommand{\R}{\mathbb{R}}

\newcommand{\PD}{\text{PD}}
\newcommand{\sub}{\sigma}

\newcommand{\hatt}{\hat{t}}

\newcommand{\modZ}{Z'}
\newcommand{\modPsi}{\Psi'}
\newcommand{\modM}{M'}

\newtheorem{thm}{Theorem}
\newtheorem{corollary}[thm]{Corollary}
\newtheorem{lemma}[thm]{Lemma}
\newtheorem{proposition}[thm]{Proposition}
\newtheorem{rem}[thm]{Remark}
\newenvironment{remark}{\bf\begin{rem}\rm}{\end{rem}} 

\title{
SINR in wireless networks and the two-parameter
  Poisson-Dirichlet process}

\author{{\bf {\bf Holger~Paul~Keeler} and  Bart{\l }omiej~B{\l }aszczyszyn}  
\thanks{INRIA-ENS, 
23 Avenue d'Italie, 75214 Paris, France}
\thanks{Bartek.Blaszczyszyn@ens.fr, keeler@ens.fr}}

\begin{document}
\maketitle
\thispagestyle{empty}

\begin{abstract}
 Stochastic geometry models of wireless networks based on Poisson
 point processes are increasingly being developed with a  focus  on
 studying various signal-to-interference-plus-noise ratio (SINR) values.
 We show that the SINR values 
experienced by a typical user with respect to different base stations
of a Poissonian cellular network are related to  a specific instance of the 
so-called two-parameter Poisson-Dirichlet process.  This process
 has many  interesting properties as well as applications in
 various fields.  We give examples of several results proved  for 
this process that are of immediate or potential interest 
in the development of  analytic tools for cellular
 networks. Some of them simplify or are akin to certain results 
that are being developed in  the network literature.
By doing this we hope to motivate further research and
use of Poisson-Dirichlet processes in this new setting.
\end{abstract}


\section{Introduction}
To derive accurate analytic tools of cellular networks, stochastic
geometry models have been developed with the almost
standard assumption that the network base stations are positioned according to a Poisson point process. The aim of these models is  often to derive distributional characteristics 
of various signal-to-interference-plus-noise ratio (SINR) values,
which are, due to information theoretic arguments, related 
to  network performance
characteristics and user quality of service metrics.
Besides tractability and `worst-case'
arguments, the Poisson assumption is justified by a recent convergence
result~\cite{hextopoi} showing that a large class of stationary
network configurations give results for functions of incoming signal strengths, such as the SINR,  as though the placement of the base stations is a Poisson process when sufficiently large log-normal shadowing is incorporated into the model~\footnote{This result is actually more general and holds under log-normal-type distributions, e.g. the Suzuki model~\cite{lognormal2014}.}.

A stochastic process known as the two-parameter Poisson-Dirichlet
process has been thoroughly studied over the years owing to
the discovery of its many interesting properties and
relations to other random structures and applications in various
fields such as population genetics, number theory, Bayesian statistics
and economics~\cite{pitman1997two,feng2010poisson}. In this letter we
detail how a specific case of this well-studied process is
equivalent  to (what we call) the
process of 
{\em  signal-to-total-interference-ratio   (STIR)} values
experienced by a typical user with respect to different base stations
of a Poissonian cellular network. The STIR process is 
trivially related to the signal-to-interference-ratio (SIR)  process
and further to the SINR one. We then list and apply results that have been derived in different settings and suggest results that may be  useful in the future.  

For related work, there is a number of Poisson-based models with a
focus on calculating the distribution of the SINR; see
~\cite{mukherjee2013analytical,blaszczyszyn2014studying} and
references therein. B\l aszczyszyn and
Keeler~\cite{blaszczyszyn2014studying}  characterized the SINR
process by obtaining its {\em factorial moment measures}. The densities
of these measure lead to the joint probability density of the order
statistics of the SINR process, which can  be used to calculate the
coverage probability under some signal combination and
interference cancellation models~\cite{blaszczyszyn2014studying}.
Invariance properties of Poisson models have been investigated in
connection to general random marks~\cite{equivalence2013} and the
special case of log-normal shadowing
marks~\cite{lognormal2014}. 
Equivalent results in
relation to the Sherrington-Kirkpatrick spin glass model
have been derived independently in physics, as detailed by
Panchenko~\cite{panchenko2013sherrington}.

The two-parameter Poisson-Dirichlet process is examined by  Pitman and
Yor~\cite{pitman1997two}, hence  it is also called the Pitman-Yor
process, though it was introduced earlier by Pearman, Pitman, and
Yor~\cite{perman1992size}. Handa~\cite{handa2009two} derived the {\em
  factorial moment density} (or {\em correlation function}) of the
process and other useful results. Kingman~\cite{kingman1992poisson}
covers the Poisson point process and its relationships to
subordinators and the original (i.e. one-parameter or Kingman's)
Poisson-Dirichlet process. In physics a related  but different 
one-parameter  process is sometimes also called the Poisson-Dirichlet
process~\cite{panchenko2013sherrington}~\footnote{It appears
as the  thermodynamic (large  system) limit in the low temperature
regime  of Derrida's  random energy model
and a key component of the so called Ruelle probability cascades,
which are used to represent the thermodynamic limit of  the
Sherrington-Kirkpatrick model for spin glasses.}, which is
  exactly our  STIR process. This process and Kingman's one are both special cases of the two-parameter Poisson-Dirichlet process.
  
 
We believe that we are the first to illustrate these connections  and in doing so it is our hope that certain results on the two-parameter Poisson-Dirichlet process will be adopted and used to develop analytic tools for studying  SINR-based characteristics in communication networks.

\section{Network model and quantities of interest}
We consider the ``typical user'' approach where one assumes a typical user is located at the origin.  On $\R^2$, we model the base stations  with a homogeneous or stationary Poisson point process  $\Phi=\{X\}$ with density $\lambda$. Define the path-loss function as
$\ell(|x|)=(K|x|)^{\beta}$,
with constants $K>0$ and $2<\beta< \infty$ assumed henceforth.   Given
$\Phi$, let $\{S_X\}_{X\in\Phi}$ be a collection of independent and
identically and \emph{arbitrarily} distributed random variables
representing  the random {\em propagation effects} (i.e fading and/or shadowing) from the origin to
$X$. Let $S$ be equal in distribution to $S_X$. In this paper we will
always (tacitly) require the moment  condition $ \E(S^{\frac{2}{\beta}})  <\infty $.

We define the propagation (loss) process~\footnote{We introduce the  propagation process for historical reasons,    in particular to be consistent with~\cite{blaszczyszyn2014studying} and papers cited therein. Otherwise the process of received powers can be considered.}, considered as a point process on the positive half-line $\mathbb{R}^{+}$, as
\begin{equation}\label{e.Y}
\Theta=\{Y\}:= \left\{\frac{\ell(|X|)}{S_X } :X\in\Phi   \right\}.
\end{equation}
\begin{lemma}\label{l.invariance}
The propagation process $\{Y\}$ is an inhomogeneous Poisson point process with intensity measure
$\Lambda_{\Theta}\left(  \left[  0,t\right)  \right)  =a t^{\frac{2}{\beta}}$,
where 
\begin{equation}\label{e.a}
a:=\frac{\lambda\pi
  \E[S^{\frac{2}{\beta}}]}{K^{2}}\,.
\end{equation}
\end{lemma}
This invariance result~\footnote{One can assume arbitrary propagation effects by setting, e.g. $S=1$, and replacing $\lambda $ with $ \lambda'=\lambda \E(S^{\frac{2}{\beta}}) $. A further generalization exists: an isotropic power-law base station density in $d$ dimensions; see, e.g.~\cite[Footnote 6]{blaszczyszyn2014studying} or~\cite[Lemma 1]{zhangdecoding}, where the propagation process is called ``path-loss process with fading''.}  has been observed a number of times; see, e.g. \cite{hextopoi} for a proof, ~\cite{blaszczyszyn2014studying} for related work, and  ~\cite{equivalence2013,lognormal2014,panchenko2013sherrington} for generalizations with random marks .

We define the {\em SINR process} on the positive half-line $\mathbb{R}^{+}$ as
\begin{equation}\label{SINR}
\Psi=\{Z\}:=\left\{\frac{Y^{-1}}{W+ (I-Y^{-1} )}:Y\in\Theta \right\},
\end{equation}
where the constant $W\geq0$ is the additive noise power, and 
\begin{equation}\label{interference}
I=\sum_{Y\in\Theta}Y^{-1},
\end{equation}
is the power received from the entire network  (so that $I-Y^{-1}$ is the interference).

To represent the {\em signal-to-total-received-power-and-noise ratio}, we define the STINR process on  $(0,1]$ as
\begin{equation}\label{MSINR}
\Psi'=\{Z'\}:=\left\{\frac{Y^{-1}}{W+ I} :Y\in\Theta \right\}.
\end{equation}
Working with $\Psi'$ is algebraically simpler and information on it gives information on  $\Psi$ by the  relation
$Z=Z'/(1- Z')$ and $Z'=Z/(1+ Z)$.

 For $ n\geq 1$, we define the factorial moment measure $\modM^{(n)}(t_1',\dots , t_n') =\modM^{(n)}((t_1',1]\times\dots \times (t_n',1])$ of the STINR process $\{\modZ \} $ as
\begin{align}
\modM^{(n)}(t_1',\dots , t_n') &=\E \left(\sum_{{(\modZ_1,\ldots,\modZ_n)\in(\modPsi)^{\times n}\atop\text{distinct}}}
\prod_{j=1}^{n}  \Ind(\modZ_{j}>t_j')\right) \label{e.Mn},
\end{align}
where $\Ind$ is an indicator function.  The equivalent measure of $\Psi=\{Z \} $ is defined by analogy but in relation to the rectangle $(t_1,\infty]\times\dots \times (t_n,\infty]$. Both measures require two integrals.  For $x\ge0$ define
\begin{equation}\label{In}
\calI_{n,\beta}(x)=\frac{2^n
\int\limits_0^{\infty} u^{2n-1}e^{-u^2-u^\beta x\Gamma(1-2/\beta)^{-\beta/2}} du
}{\beta^{n-1}(C'(\beta))^n(n-1)!}
\end{equation}
where 
 \begin{equation}\label{e.C}
 C'(\beta):=2\pi/(\beta\sin(2\pi/\beta))=\Gamma(1+2/\beta)\Gamma(1-2/\beta),
 \end{equation}
  and $\Gamma$ is the gamma function. Note that   $\calI_{n,\beta}(0)=2^{n-1}/(\beta^{n-1}[C'(\beta)]^n)$.  
For all $x_i\ge0$ define
\begin{align}\label{Jn}
& \calJ_{n,\beta}(x_1,\dots,x_n)= \frac{(1+\sum_{j=1}^{n}   x_j) }{n} \nonumber \\ 
&\times  \int\limits_{[0,1]^{n-1}}  \frac{  \prod_{i=1}^{n-1}   v_i^{i(2/\beta+1)-1}(1-v_i)^{2/\beta}  }   {  \prod_{i=1}^n  (x_i+\eta_i)}
   dv_1\dots
dv_{n-1},
\end{align}
where 
\begin{align}
\eta_1&= v_1v_2\dots v_{n-1}, \quad \eta_2= (1-v_1)v_2\dots v_{n-1},\nonumber \\
\eta_3&= (1-v_2)v_3\dots v_{n-1},\quad\cdots,\quad\eta_n= 1- v_{n-1}.
\end{align}
For reasonably low $n$ (i.e. $n\leq 20$), both these integrals are numerically tractable~\cite{paul_matlab_moments}.
Let $\hatt_{i}=\hatt_{i}(t_1',\dots,t_n'):=t_i'/(1-\sum\limits_{j=1}^n t_j')$  and the define the $n$-dimensional unit simple
\[
\Delta_n=\{(t_1',\dots t_n'): t_1', \dots  t_n'\ge 0,\; t_1'+ \dots + t_n'  \leq 1\},
\]
and $\Ind_{\Delta_n}$ denotes the corresponding indicator function. 
 We now present the factorial moment measures~\cite{blaszczyszyn2014studying}.
\begin{proposition}\label{mainResult}
For $t_i'\in(0,1]$, the factorial moment measure of order $n \ge1$ of the STINR process  (\ref{MSINR})  satisfies
\begin{align} \label{momMeasure}
&\modM^{(n)}(t_1',\dots  t_n') = n!  \left( \prod\limits_{i=1}^{n}\hatt_i^{-2/\beta} \right) \nonumber \\
&\times \calI_{n,\beta}(Wa^{-\beta/2}) 
 \calJ_{n,\beta}(\hatt_1,\dots,\hatt_n) \Ind_{\Delta_n} (t_1',\dots,t_n')  .
\end{align}
Furthermore, for $t_i\in(0,\infty)$ the SINR process  (\ref{SINR}) has the moment measure
\begin{equation}\label{e.MmodM}
M^{(n)}\left(t_1,\dots t_n,\right)  = \modM^{(n)}\left(t_1',\dots, t_n'\right),
 \end{equation}
 where $t_i=t'_i/(1- t'_i)$ and $t'_i=t_i/(1+ t_i)$.

\end{proposition}

Let $ \modM^{(n)}_0 $ and  $ M^{(n)}_0$  respectively denote the
factorial moment measures of the STIR and SIR processes, i.e.,  $\modM^{(n)} $ and  $ M^{(n)}$ with   $W=0$.  Hence
\begin{align}
\modM^{(n)}(\cdot)&= \bar{ \calI}_{n,\beta}(Wa^{-\beta/2})  \modM^{(n)}_0(\cdot)\label{noiseM},\\
M^{(n)}(\cdot)&= \bar{ \calI}_{n,\beta}(Wa^{-\beta/2})M^{(n)}_0(\cdot) ,
 \end{align}
where
\begin{equation}
\bar{ \calI}_{n,\beta}(Wa^{-\beta/2})=\frac{ \calI_{n,\beta}(Wa^{-\beta/2})} { \calI_{n,\beta}(0)} .
 \end{equation}
This ability to separate the noise term in the factorial moment measures (and densities) is  convenient and is reminiscent of factoring out the noise term in the distribution of the SINR under Rayleigh fading, an assumption that is not required, however, in our present setting.

\section{Two-parameter Poisson-Dirichlet process}
One way to define the two-parameter Poisson-Dirichlet process~\cite{pitman1997two} for two given parameters $0\leq \alpha < 1$ and $\theta>-\alpha$ is to first introduce a sequence of random variables $\{\tilde{V}_i\}$ by
\begin{equation}\label{e.tildeV}
\tilde{V}_1=U_1 ,\quad \tilde{V}_i=(1-U_1)\dots (1-U_{i-1})U_i , \quad i\geq2,
 \end{equation}
where $U_1, U_2,\ldots$ are independent beta variables such that each $U_i$
has $B(1-\alpha, \theta+i\alpha)$ distribution. Note that 
$\sum_{i=1}^{\infty}\tilde V_i=1$ with probability
one.  Denote the decreasing order statistics $\{\tilde{V}_{(i)}\}$ of $\{\tilde{V}_i\}$ by
$\{V_i\}$ ($V_1\geq V_2\geq \dots$), then define the two-parameter
Poisson-Dirichlet distribution with parameters $\alpha$ and $\theta$,
abbreviated as $\PD(\alpha,\theta)$, to be the distribution of
$\{V_i\}$. By considering $\{V_i\}$ (or equivalenently
$\{\tilde{V}_i\}$) as atoms of a point process, we see
$\PD(\alpha,\theta)$ as a distribution of a point process. The above
approach of defining the $\PD(\alpha,\theta)$ distribution is related
to problems on so-called size-biased sampling and stick-breaking or
the residual allocation model, where $\PD(\alpha,\theta)$  plays a
central role. In fact, the distribution of $\{\tilde{V}_i\}$ coincides with that of  
the size-biased permutation of $\{V_i\}$~\cite{pitman1996under}. 

Another way to define a Poisson-Dirichlet
process~\cite{pitman1997two}, more aligned with our setting, is to use
the concept of a subordinator 
having almost surely increasing trajectories. 
 For $s\geq 0$, let $\sub_s$ be a subordinator, and, assuming it has zero drift, its Laplace transform is
\begin{equation}\label{laplacesub}
\E[\exp(-z \sub_s)]=\exp\left[-s\int_0^{\infty}  (1-e^{-z r}) \Lambda (dr)\right],
 \end{equation}
where $\Lambda$ is a measure on  $(0,\infty)$,
called the L\'evy measure, characterizing the subordinator without drift. 
When multiplied by $s$, this measure ($s\Lambda(dr)$)   
can be identified with the intensity measure of
  the Poisson point process of jumps the subordinator makes
 in the interval $(0,s)$. Let us order and denote these jumps by
 $V_1(\sigma_s)\geq V_2(\sigma_s)\geq \dots$. Clearly $\sub_s=\sum_{i=1}^{\infty}
 V_i(\sigma_s) $.

Let $0<\alpha < 1$, and then $\sub_s$ is called an $\alpha$-stable subordinator if
$\Lambda(dr)=D r^{-\alpha-1}\,dr$ for 
some constant $D>0$,  which implies $\E[\exp(-z \sub_s)]=\exp[- s D
   \Gamma (1-\alpha) z^{\alpha}]$. {A crucial
observation~\cite[Proposition 6]{pitman1997two} says that  for any
$s>0$ the sequence
$\{V_1(\sub_s)/\sub_s,V_2(\sub_s)/\sub_s,\dots\}$ has
$\PD(\alpha,0)$ distribution~\footnote{To define the original
   Kingman's Poisson-Dirichlet process with  $\PD(0,\theta)$
   distribution, a subordinator known as the
   Moran~\cite{kingman1992poisson} or Gamma
   subordinator~\cite{pitman1997two}  is used. A subtle combination of
   this and the $\alpha$-stable subordinator is used to define the
   $\PD(\alpha,\theta)$ process, but we omit  the details for
   brevity.}.

Set $s=1$ and
the constants 
$\alpha=2/\beta$ and $D=a$, 
where $a$ is given by  (\ref{e.a}).
Then we see that 
the jumps of the subordinator in the interval $(0,s)$ can be identified with the power
values of the signals from all the base stations or, equivalently, the
inverse values of the propagation process $\Theta$
which, in view of Lemma~\ref{l.invariance}, is
an inhomogeneous Poisson process with intensity measure
$(2a/\beta) t^{-1-2/\beta}dt$.
Consequently, $\sigma_1$ represents the interference in  our Poisson
network model, and its  Laplace  transform is 
$\E[\exp(-z \sub_1)]=\E[\exp(-z I)]=\exp[-  a\Gamma (1-2/\beta)
  z^{2/\beta}]$.

In other words,  the subordinator representation of the Poisson-Dirichlet process  $\PD(\alpha,0)$ 
(\cite[Proposition 6]{pitman1997two}) relates this process
to our STIR process. More precisely, 
denote the {\em increasing} orders statistics of  $\{ Y\}$ by $\{
Y_{(i)}\}$, such that $Y_{(1)}\leq Y_{(2)} \leq \dots$,  and the {\em
  decreasing} order statistics of $\{ Z'\}$ by $\{ Z_{(i)}'\}$. Then
we have the following relation, which is a key observation of this letter.
\begin{proposition} Assume $W=0$. Then the sequence $\{Z'_{(i)}\}$ is
  equal in distribution to $\{V_i\}$ for $\alpha=2/\beta$ and $\theta=0$. In other words, the STIR  process $\Psi'$ is a $\PD(2/\beta,0)$ point process.
\end{proposition}
The fact that $\{ \tilde{V}_i\}$, defined
  in~(\ref{e.tildeV}), 
form a size-biased permutation of $\{V_i\}$ can be interpreted as follows regarding our STIR process. 
\begin{remark}
Assume  when the typical user is choosing its serving base station that, instead of looking for the strongest received signal $Y^{-1}_{(1)}$,  it makes a randomized decision, picking a base station $i$ with a bias proportional to  $Y^{-1}_{i}$ (hence stronger stations have more chance to be selected). 
Then its STIR, with respect to the chosen station, has the
distribution of $\tilde V_1$, i.e. $B(1-2/\beta,2/\beta)$. Suppose now
that another user positioned with the typical one and subject to the
same propagation effects makes its choice of the serving base station
by applying the same randomized procedure but excluding the station
already selected by the first user. Then the join distribution of the
STIR's experienced by these two users is equal to that of the random
pair $(\tilde V_1,\tilde V_2)$, which can be easily derived
from~(\ref{e.tildeV}). The above {\em randomized access policy}, and the corresponding evaluation of the STIR values, which can be extended to an arbitrary number of users, is of potential interest for managing user hotspots.
\end{remark}

\section{Some useful results}
 Appropriately adapted for this setting, we list some interesting results of the $\PD(2/\beta,0)$ distribution.  The first result~\cite[Proposition 8]{pitman1997two} applied here shows that the ratio of successive STINR values have beta distributions. 
\begin{proposition}
For the STINR process $\Psi'$ ($W\geq 0$), the random variables
\begin{equation}
R_i:= \frac{Z'_{(i+1)}}{Z'_{(i)}}=\frac{Y_{(i)}}{Y_{(i+1)}} 
\end{equation}
have, respectively, $B( 2i/\beta,1)$ distributions such that $\Prob( R_i \leq r)=r^{i 2/\beta}$ (for $0\leq r \leq 1$). Moreover, $\{R_i\}$ are mutually independent. 
\end{proposition}
The fact that each $R_i$ is a ratio of $Y_{(i)}$ values indicates that
this  result (proved in~\cite{pitman1997two} under assumption $W=0$)
is invariant of the noise term $W$. This applies also to the next
result, which involves the following variables
\begin{align}
A_i:&=\frac{Z'_{(1)}+\dots +  Z'_{(i)}}{Z'_{(i+1)}} = \frac{ Y_{(1)}^{-1}+\cdots + Y_{(i)}^{-1}}{Y_{(i+1)}^{-1}}.  \\
\Sigma_i:&=\frac{Z'_{(i+1)}+Z_{i+2}'+\dots}{Z_i'}= \frac{ Y_{(i+1)}^{-1}+ Y_{(i+2)}^{-1}+\cdots  }{Y_{(i)}^{-1}}
\end{align}
defined for $i=1, 2, \dots$.
For $\gamma\geq 0$ let 
\begin{align}
\phi_{\beta}(\gamma)&:=\frac{2}{\beta}\int_1^{\infty} e^{-\gamma x} x^{-2/\beta -1 } dx , \\
\psi_{\beta}(\gamma)&:=\Gamma(1-2/\beta) \gamma^{2/\beta} + \phi_{\beta}(\gamma) .
\end{align}
The next proposition follows~\cite[Proposition 11]{pitman1997two}.
\begin{proposition}\label{ASigma}
Consider the STINR process $\Psi'$  ($W\geq 0)$. 
Then 
$1/Z'_{(i)}=1+A_{i-1} + \Sigma_i$
where $A_{i-1}$ is distributed as the sum of $i-1$ independent copies of $A_1$, with the characteristic function $\E[e^{-\gamma A_{i-1}}]=(\phi_{\beta}(\gamma))^{i-1}$; $\Sigma_{i}$ is distributed as the sum of $i$ independent copies of $\Sigma_1$, with the characteristic function $\E[e^{-\gamma \Sigma_{i}}]=(\psi_{\beta}(\gamma))^{-i}$; and $A_{i-1}$ and $\Sigma_i$ are independent. 
\end{proposition}
\begin{remark}
By observing that $\Sigma_i^{-1}=Y_{(i)}^{-1}/(Y_{(i+1)}^{-1}+Y_{(i+2)}^{-1}+\ldots) $,
the above result in the setting of {\em successive-interference cancellation} (with no noise, $W=0$) can be compared to a result ~\cite[Theorem 1]{zhangdecoding} and its generalization~\cite[Proposition 21]{blaszczyszyn2014studying} on the ratio of the $k\,$th strongest propagation process and a successively reduced interference term. Moreover, the ratio of independent random variables  $(1+A_{i-1})/\Sigma_i = (Y_{(1)}^{-1}+\cdots + Y_{(i)}^{-1})/(Y_{(i+1)}^{-1}+Y_{(i+2)}^{-1}+\cdots) $ relates the above result to a recent {\em signal combination} model in the STIR ($W=0$) scenario~\cite{blaszczyszyn2014studying}. The difference between the STINR and STIR results suggests  that the  noise term $W$ can add a significant layer of complexity to the models.
\end{remark}
Proposition \ref{ASigma} leads to a Laplace transform result (cf~\cite[Corollary 12]{pitman1997two}), which
can be compared to a previous observation~\cite[Remark
  18]{blaszczyszyn2014studying}.
\begin{corollary}
The inverse of the $i\,$th strongest STIR ($W=0$) value,  $1/Z'_{(i)}$, has the Laplace transform
\begin{equation}
\E[e^{-\gamma /Z'_{(i)} }]=e^{-\gamma} (\phi_{\beta}(\gamma))^{i-1} (\psi_{\beta}(\gamma))^{-i}.
\end{equation}
\end{corollary}
Furthermore, a previous result~\cite[Corollary 7]{kcovsingle} gives an expression for the tail of the distribution function of   the $i\,$th strongest STINR ($W\ge 0$) value.

Remarkably, the interference $I$, as defined by (\ref{interference}),
and noise $W$ can be recovered from the STINR processes. Indeed, the
first statement of the following result is trivial, while the second
one can be proved using the same arguments as~\cite[Proposition 10]{pitman1997two}.
\begin{proposition}
For the STINR process ($W\geq 0$), 
$W/I= \Bigl(\sum_{i=1}^{\infty}Z'_{(i)}\Bigr)^{-1} -1$,
and $W+I= (L/a)^{-\beta/2}$,  where the limit
$L:=\lim_{i \rightarrow \infty} i( Z'_{(i)})^{2/\beta}$, both exists almost surely and for all $p$-means with $p \geq 1$. 

\end{proposition}

The densities of the factorial moment measures $M'^{(n)}$ of the STINR process 
can be used to find an expression for the joint probability density of
the order statistics of the STINR process~\cite[Proposition
  20]{blaszczyszyn2014studying}. But, the expression on the
right-hand-side of~(\ref{momMeasure}) appears too unwieldy to
differentiate with respect to more than a couple variables. However,
using the representation~(\ref{e.tildeV}), the factorial moment density of the
$\PD(\alpha,\theta)$ process was derived in closed-form~\cite[Theorem
  2.1]{handa2009two}, which implies the following new result for our STINR process.
For $n\ge0$ denote $c_{n,\alpha,\theta}=\prod_{i=1}^n\Gamma(\theta
+1+(i-1)\alpha)/(\Gamma((1-\alpha)\Gamma(\theta +i\alpha))$;
in particular  $c_{n,2/\beta,0}=   (2/\beta)^{n-1} \Gamma(n )/(\Gamma(2n/\beta )\Gamma(1-2/\beta )^n)$.
\begin{proposition}\label{newmu_n}
For the STINR process $\Psi'$ $(W\geq 0)$, the $n\,$th factorial
moment density is given by
\begin{align} 
&\mu^{(n)}(t_1',\dots  t_n') :=(-1)^n
\frac{\partial^{n}\modM^{(n)}(t_1',\dots  t_n')}{\partial
  t_1'\dots\partial t_n' }
\label{mu_n}
\\
&\hspace{-1em}= c_{n,2/\beta,0}\,\bar{ \calI}_{n,\beta}(Wa^{-\beta/2})
 \Bigl( \prod\limits_{i=1}^{n}t_i'^{-(2/\beta+1)} \Bigr)
\Bigl(1- \sum\limits_{j=1}^{n}t_j' \Bigr)^{2/\beta n-1}\hspace{-3em}\,
\nonumber
\end{align}
for $(t_1',\dots,t_n')$ in ${\Delta_n}$ and 0 otherwise.

\end{proposition}
This result follows from including the noise term $W $, via $(\ref{noiseM})$, and using ~\cite[Theorem 2.1]{handa2009two}.  We are unaware of anybody
    showing  the equivalence of Propositions \ref{mainResult} and  \ref{newmu_n}, either by
    differentiating the    measure (\ref{momMeasure}) or integrating the density (\ref{mu_n}). 

 Another way for calculating the joint density of the order statistics
 of the STIR process is offered by the result~\cite[Theorem 5.4]{handa2009two}, where the two-parameter Dickman function  was introduced as
\begin{equation}\label{dickman}
\rho_{\alpha,\theta}(s):= \Prob(V_1<1/s)\,
\end{equation}
 where $V_{1}$ is the largest value of the $\PD(\alpha,\theta)$
 process, 
which can be computed as follows:
\begin{equation}\label{e.rhoat}
\rho_{\alpha,\theta}(s):= \sum\limits_{n=0}^{\infty} \frac{(-1)^n
  c_{n,\alpha,\theta} }{n!} I_{n,\alpha,\theta }(s),
\end{equation}
\vspace{-1.5ex}
where 
\vspace{-1.2ex}
\begin{equation}\label{e.Inat}
I_{n,\alpha,\theta}(s)= \int_{\Delta_n}  \prod_{i=1}^{n} \frac{\Ind_{[1,\infty)}(s t_i')}{t_i'^{\alpha+1}}   \Bigl(1- \sum\limits_{j=1}^{n}t_j'  \Bigr)^{\theta+\alpha n-1}\hspace{-2em} dt_1'\dots dt_n',
\end{equation}
and for $n=1,2, \dots$,
with $I_{n,\alpha,\theta}(s) = 0$ whenever $n>s$, which makes that the
right-hand-side of~(\ref{e.rhoat}) is actually a finite sum;
cf~\cite[Section 4]{handa2009two}.
\begin{remark}
Recall that when $\alpha=2/\beta$ and $\theta=0$, $V_1$ is equal in
distribution to the strongest STIR value $Z'_{(1)}$ and
thus (\ref{e.rhoat}) should be
compared to~\cite[Corollary 7, with $k=1$, valid for $W\ge0$]{kcovsingle}
or~\cite[Theorem 1, taken for single-tier network, valid for SINR values greater
than one]{dhillon2012modeling}.
\end{remark}
\begin{proposition}
For the STIR process $\Psi'$ $(W=0)$ and for each $m=1,2,\dots$, the joint probability density of $\{Z'_{(1)},\dots, Z'_{(m)}  \}$ is given by
\begin{align} 
f_{m,\beta}&(t_1',\dots  t_m') :=\, c_{m,2/\beta,0} \Bigl(
\prod\limits_{i=1}^{m}t_i'^{-(2/\beta+1)} \Bigr)
\! \Bigl(1- \sum\limits_{j=1}^{m}t_j'  \Bigr)^{2m/\beta-1}  \nonumber\\
&\times  \rho_{2/\beta,2m/\beta}\Bigl(\frac{1- \sum_{j=1}^{m}t_j'}{t'_m} \Bigr)\Ind_{\Delta_m}(t_1',\dots,t_m')
\Ind_{\{t_1'> \dots >t_n'\}}\,.
\end{align}
\end{proposition}
The two non-zero parameters of the Dickman function explains why $\PD(\alpha,\theta)$  is needed in~(\ref{dickman}), and not just $\PD(\alpha,0)$.  

\section{Conclusion}
We showed the relationship between the SINR process, which is an
important object in the  study of  the performance of cellular
networks, and the  two-parameter Poisson-Dirichlet process. We
presented some results recently proved for the former process, which
have interesting interpretation in terms of the  STINR process (easily
related to the SINR one). Our goal is to encourage  further research aimed at building bridges
between these two, until now, separate research areas.

\paragraph*{Acknowledgements}
We thank M.K. Karray for carefully proofreading the manuscript.

{\small
\pdfbookmark[0]{References}{References} 
\bibliographystyle{IEEEtran}

}
\end{document}